\RequirePackage{fix-cm}
\documentclass[smallextended]{svjour3}       % onecolumn (second format)
\smartqed  % flush right qed marks, e.g. at end of proof

\usepackage{array}
\usepackage{amssymb}
\usepackage{booktabs}
\usepackage{subfigure}
\usepackage{graphicx}
\begin{document}

\title {Two efficient computational algorithms to solve the singularly perturbed Lane-Emden problem}


\titlerunning{The singularly perturbed Lane-Emden problem}        % if too long for running head

\author{Kourosh Parand$^{1}$ \and Amin Ghaderi$^{2}$ }


\institute{\at $^1$ Department of Computer Sciences, Shahid Beheshti University, G.C., Tehran, Iran, \\  Department of Cognitive Modelling, Institute for Cognitive and Brain Sciences, Shahid Beheshti University, Tehran, Iran. Email: k\_parand@sbu.ac.ir.
\\
\at {$^2$} Department of Cognitive Modelling, Institute for Cognitive and Brain Sciences, Shahid Beheshti University, Tehran, Iran,  Email: amin.g.ghaderi@gmail.com}
\date{Received: date / Accepted: date}

\maketitle

\begin{abstract}
In this paper, we decide to compare two new approaches based on Rational and Exponential Bessel functions (RBs and EBs) to solve several well-known class of Lane-Emden type models. The problems, which define in some models of non-Newtonian fluid mechanics and mathematical physics, are nonlinear ordinary differential equations of second-order over the semi-infinite interval and have singularity at $x=0$. We have converted the non-linear Lane-Emden equation to  a sequence of  linear equations by utilizing the quasilinearization method (QLM) and then, these linear equations have been solved  by RBs and EBs collocation-spectral methods. Afterward, the obtained results are compared with the solution of other methods for demonstrating the efficiency and applicability of the proposed methods.
\end{abstract}

\section{Introduction}\label{sec1}
The investigation of a singular initial/boundary value non-linear differential equations of second-order have been attracted by some astrophysicist, mathematicians and physicists. Lane-Emden type equations describe the temperature variation of a spherical gas cloud under the mutual attraction of its molecules and subject to the laws of classical thermodynamics. 
Let $P(r)$ denote the total pressure at a distance $r$ from the center of spherical gas cloud. The total pressure is due to the usual gas pressure and a contribution from radiation:
\begin{eqnarray}
\nonumber P=(\frac{1}{3}\xi T^{4}+\frac{RT}{v}),
\end{eqnarray}
where $\xi, T , R$ and υ are respectively the radiation constant, the absolute temperature, the gas constant, and the specific volume, respectively\cite{refff03,refff04}. Let $M(r)$ be the mass within a sphere of radius $r$ and $G$ the constant of gravitation. The equilibrium equation for the configuration are
\begin{eqnarray}\label{lbl3}
&&\frac{dP}{dr}=-\rho \frac{GM(r)}{r^{2}},\\
&&\nonumber \frac{dM(r)}{dr}=4\pi \rho r^{2},
\end{eqnarray}
where $\rho$, is the density, at a distance $r$, from the center of a spherical star. Eliminating $M$ yields of these equations results in the following equation, which as should be noted, is an equivalent form of the Poisson equation \cite{refff04,refff05}:
\begin{eqnarray}
\nonumber\frac{1}{r^2}\frac{d}{dr}(\frac{r^2}{\rho}\frac{dP}{dr})=-4\pi G\rho.
\end{eqnarray}
We already know that in the case of a degenerate electron gas, the pressure and density are$ \rho = P^{\frac{3}{5}}$, assuming that such a relation exists in other states of the star, we are led to consider a relation of the form $P = K\rho^{1+\frac{1}{m}}$, where $K$ and $m$ are constants.

We can insert this relation into Eq.(\ref{lbl3}) for the hydrostatic equilibrium condition and, from this, we can rewrite the equation as follows:
\begin{eqnarray}
\nonumber\Bigg[\frac{K(m+1)}{4\pi G}\lambda^{ \frac{1}{m}-1}\Bigg]\frac{1}{r^2}\frac{d}{dr}(r^{2}\frac{dy}{dr})=-y^{m},
\end{eqnarray}
where $\lambda$ represents the central density of the star and y denotes the dimensionless quantity, which are both related to $\rho$ through the following relation \cite{refff01,refff05}:
\begin{eqnarray}
\nonumber\rho=\lambda y^{m}(x),
\end{eqnarray}
and let
\begin{eqnarray}
&&\nonumber r=ax,\\
&&\nonumber a=\Bigg[\frac{K(m+1)}{4\pi G}\lambda^{\frac{1}{m}-1}\Bigg]^{\frac{1}{2}}.
\end{eqnarray}
Inserting these relations into our previous relation we obtain the Lane-Emden equation \cite{refff04,refff05}:
\begin{eqnarray}
\nonumber\frac{1}{x^2}\frac{d}{dx}(x^2\frac{dy}{dx})=-y^{m},
\end{eqnarray}
now, we will have the standard Lane-Emden equation with $f(x, y) = y^{m}$
\begin{eqnarray}\label{lbl4}
y''(x)+\frac{2}{x}y'(x)+y^{m}(x)=0,~~x>0,
\end{eqnarray}
the initial conditions are as follows
\begin{eqnarray}\label{lbl6}
y(0)=1~,~y'(0)=0.
\end{eqnarray}
The values of $m$, which are physically interesting, lie in the interval [0, 5]. The main difficulty in analyzing this type of equation is the singularity behaviour occurring at $x = 0$.

As it has been mentioned in the literature review, the solutions of the Lane-Emden equation could be exact only for $m = 0, 1$ and $5$. For the other values of $m$, the Lane- Emden equation is to be integrated numerically \cite{refff05}. Thus, we decided to present a new and efficient technique to solve it numerically for $m = 0.5, 1.5, 2, 2.5, 3, 3.5 $ and $4$.

\subsection{Previous  works}
Recently, many analytical, semi- analyticaland and numerical techniques have been applied to solve Lane-Emden equations. The main difficulty arises in the singularity of the equations at $x = 0$. We have introduced several techniques as follow:

Bender et al. \cite{refff06} proposed a new perturbation technique based on an artificial parameter $\delta$, the method is often called $\delta$-method.  Wazwaz \cite{refff08} employed the Adomian decomposition technique with an alternate framework designed, J.H. He \cite{refff11} employed Ritz’s method to obtain an analytical solution, Parand et. al. \cite{refff15,refff16,refff17,refff18,refff19} applied pseudo-spectral method based on rational Legendre functions, Sinc collocation method, the Lagrangian method based on modified generalized Laguerre function, Hermite function collocation method and meshless collocation method based on Radial basis function (RBs) as numerical solution, Ramos \cite{refff20,refff21,refff22,refff23} presented linearzation methods to utilize an analytical solutions and globally smooth solutions, developed piecewise-adaptive decomposition methods, obtained series solutions of the Lane-Emden type equation, Yousefi \cite{refff24}  applied Legendre Wavelet approximations and used integral operator and converted Lane-Emden equations to integral equations, Chowdhury and Hashim \cite{refff25} used analytical solutions of the generalized Emden- Fowler type equations by homotopy perturbation method (HPM), Aslanov  \cite{refff27} introduced a further development in the Adomian decomposition technique, Dehghan and Shakeri \cite{refff28} investigated Lane-Emden equations by applying the variational iteration method, Marzban et al. \cite{refff30} used a method based upon hybrid of block-pulse functions and Lagrange interpolating polynomials together with the operational integration matrix to approximate solution of the problem, Adibi and Rismani in \cite{refff31}  proposed the approximate solutions of singular the Lane-Emden via modified Legenre-spectral method, Karimi vanani and Aminataei \cite{refff32}  provided a numerical method which produces an approximate polynomial solution, thesy used an integral operator and convert Lane-Emden equations into integral equations then, convert the acquired integral equations into a power series and finally, transforming the power series into pad\'{e} series form, Kaur et al.  \cite{refff33} obtained the Haar wavelet approximate solution. 

So, the other researchers trying to solving the Lane-Emden type equations with several methods, For example, A Yildirım and \"{O}zi\c{s} \cite{refff34,refff35} by using HPM and VIM methods, Benko et al. \cite{refff36} by using Nystr\"{o}m method, Iqbal and javad \cite{refff37}  by using Optimal HAM, Boubaker and Van Gorder \cite{refff38} by using boubaker polynomials expansion scheme, Da\c{s}cıo\v{g}lu and Yaslan \cite{refff39} by using Chebyshev collocation method, Y\"{u}zba\c{s}ı \cite{refff40,refff41} by using Bessel matrix and improved Bessel collocation method, Boyd \cite{{refff42}} by using Chebyshev spectral method, Bharwy and Alofi \cite{refff43} by using Jacobi-Gauss collocation method, Pandey et al. \cite{refff44,refff45} by using Legendre and Brenstein operation matrix, Rismani and monfared \cite{refff46} by using Modified Legendre spectral method, Nazari-Golshan et al. \cite{refff47} by using Homotopy perturbation with Fourier transform, Doha et al. \cite{refff48} by using second kind Chebyshev operation matrix algorithm, Carunto and bota \cite{refff49} by using Squared reminder minimization method, Mall and Chakaraverty \cite{refff50} by using Chebyshev Neural Network based model, G\"{u}rb\"{u}z and sezer \cite{refff51} by using Laguerre polynomial and Kazemi-Nasab et al. \cite{refff52} by using Chebyshev wavelet finite difference method. In this paper, we attempt to introduce a new method, based on RBF-DQ for solving non-linear ODEs.

The rest of this paper is arranged as follows:  

Section \ref{sec2}  introduces new rational and exponential Bessel functions (RBs and EBs) and their properties. Section \ref{sec3} describes a brief formulation of quasilinearization method (QLM) introduced by \cite{QLM03}. In section \ref{sec4} at first, by utilizing QLM over Lane-Emden equation  a sequence of linear differential equations is obtained, then at each iteration solve the linear differential equation by RBs and EBs collocation methods that we name RBs-QLM and EBs-QLM methods. Comparison between these two methods  with some well-known results in section \ref{sec5}, show that using rational functions is highly accurate, and we also  describe our results via tables and figures. Finally, we give a brief conclusion in section \ref{sec6}.

\section{Properties of Rational and Exponential Bessel Functions}\label{sec2}
The Bessel functions arise in many problems in physics possessing cylindrical symmetry, such as the vibrations of circular drumheads and the radial modes in optical fibers. Bessel functions are usually defined as a particular solution of a linear differential equation of the second order which known as Bessel's equation. Bessel functions first defined by the Daniel Bernoulli on heavy chains (1738) and then generalized by Friedrich Bessel. More general Bessel functions were studied by Leonhard Euler in (1781) and in his study of the vibrating membrane in (1764) \cite{reffff01,reffff02}.
%\subsection{Bessel polynomials}

Bessel differential equation of order $n\in\mathbb{R}$ is:
\begin{eqnarray}\label{lbl7}
x^2\frac{d^2y(x)}{dx^2}+x\frac{dy(x)}{dx}+(x^2-n^2) y(x)=0,~~x\in(-\infty,\infty).
\end{eqnarray}

One of the solutions of equation (\ref{lbl7}) by applying the method of Frobenius as follows \cite{reffff03}:
\begin{eqnarray}\label{lbl8}
J_{n}(x)=\sum_{r=0}^{ \infty}\frac{(-1)^r}{r!(n+r)!}(\frac{x}{2})^{2r+n},
\end{eqnarray}
where series (\ref{lbl8})  is convergent for all $x\in(-\infty,\infty)$.

polynomials has been introduced as follows \cite{reffff11,reffff12}:
\begin{eqnarray}\label{lbl9}
B_{n}(x)=\sum_{r=0}^{[\frac{N-n}{2}]}\frac{(-1)^r}{r!(n+r)!}(\frac{x}{2})^{2r+n},~~x\in[0,1].
\end{eqnarray}
where $n\in\mathbb{N}$, and $N$ is the number of basis of Bessel polynomials.

\subsection{Rational Bessel Functions}
The new basis functions, "Rational Bessel functions  (RBs)" denote by $RB_{n}(x,L)$ which  are generated from well known Bessel polynomials by using the algebraic mapping $\phi(x)=\frac {x}{x + L}$, as follow:
\begin{eqnarray}
\nonumber RB_{n}(x,L)=B_{n}(\phi(x))~,~~n=0, 1,\cdots,N\\ 
\nonumber
\end{eqnarray}
or
\begin{eqnarray}\label{lbl10}
RB_{n}(x,L)=\sum_{r=0}^{[\frac{N-n}{2}]}\frac{(-1)^r}{r!(n+r)!}(\frac{x}{2(x + L)})^{2r+n},~~n=0, 1,\cdots,N
\end{eqnarray}
where $x\in [0,\infty)$, $B_{n}(x)$ is Bessel polynomials of order $n$, and the constant parameter $L>0$ is a scaling/stretching factor which can be used to fine tune the spacing of collocation points. For a problem whose solution decays at infinity, there is an effective interval outside of which the solution is negligible, and collocation points which fall outside of this interval are essentially wasted. On the other hand, if the solution is still far from negligible at the collocation points with largest magnitude, one cannot expect a very good approximation. Hence, the performance of spectral methods in unbounded domains can be significantly enhanced by choosing a proper scaling parameter such that the extreme collocation points are at or close to the endpoints of the effective interval \cite{reffff17}. Boyd \cite{refdd26} offered guidelines for optimizing the map parameter $L$ for rational Chebyshev functions, which is also useful for RBs.

Let us define $\Gamma=\{x|~0\leq x < \infty  \}$ and \\
$L^{2}_{w_{r}}(\Gamma)=\{~v :\Gamma \rightarrow \mathbb{R}| v$ is measurable and $\parallel v \parallel_{w_{r}} < \infty \}$, where
$$\parallel v \parallel_{w_{r}}=\left(\int^{\infty }_{0}|v(x)|^{2}w_{r}(x,L)dx\right)^{1/2},$$
with $w_{r}(x,L)=\frac{L}{(x+L)^{2}}$, is the norm induced by inner product of the space $L^{2}_{w_{r}}(\Gamma)$ as follows:
$$\langle v,u\rangle_{w_{r}}=\int^{\infty }_{0}{v(x)u(x)w_{r}(x,L)}dx.$$

Now, suppose that\\
$$\mathfrak{S}=~span\{RB_{0}(x), RB_{1}(x),\dots, RB_{N}(x)\},$$
where $\mathfrak{S}$ is a finite-dimensional subspace of $L^{2}_{w}(\Gamma)$,  $dim( \mathfrak{S}) = N+1$, so $\mathfrak{S}$ is a closed subspace of $L^{2}(\Gamma)$. Therefore, $\mathfrak{S}$ is a complete subspace of $L^{2}(\Gamma)$. Assume that $f(x)$ be an arbitrary element in $L^{2}(\Gamma)$. Thus $f(x)$ has a unique best approximation in $\mathfrak{S}$ subspace, say $\hat{b}(x)\in \mathfrak{S}$, that is
\begin{eqnarray}
\nonumber\exists~ \hat{b}(x)\in\mathfrak{S}, ~~~ \forall ~b(x)\in \mathfrak{S}(x),~~\parallel f(x)-\hat{b}(x)\parallel_{w_{r}} \leq \parallel f(x)-b(x)\parallel_{w_{r}}.
\end{eqnarray}
Notice that we can write $b(x)$ vector as a combination of the basis vectors of $\mathfrak{S}$ subspace. 

We know function of $f(x)$  can be expanded by $N+1$ terms of RB as: 
\begin{eqnarray}
f(x)=f_{N}(x)+R(x),
\nonumber\end{eqnarray}
that is
\begin{equation}\label{lbl11}
f_{N}(x)=\sum^{N}_{n=0}{a_{n}RB_{n}(x)}=A^{T}RB(x),
\end{equation}
where  $RB(x)$ is vector $[RB_{0}(x), RB_{1}(x),\cdots, RB_{N}(x)]^{T}$ and $R\in\mathfrak{S}^{\perp}$ that ${\mathfrak S}^{\perp}$ is the orthogonal complement. So  $f(x)-f_{N}(x)\in\mathfrak{S}^{\perp}$ and $b(x)\in\mathfrak{S}$ are orthogonal which we denote it by:
\begin{eqnarray}
\nonumber f(x)-f_{N}(x)\perp b,
\end{eqnarray}
thus $f(x)-f_{N}(x)$ vector is orthogonal over all of basis vectors of $\mathfrak{S}$ subspace as:
\begin{eqnarray}
\nonumber \langle f(x)-f_{N}(x),RB_{i}(x)\rangle_{w_{r}}=\langle f(x)-A^{T}RB(x),RB_{i}(x)\rangle_{w_{r}}=0,~i=0, 1,\cdots, N,
\end{eqnarray}
hence
\begin{eqnarray}
&&\nonumber\langle f(x)-A^{T}RB(x),RB^{T}(x)\rangle_{w_{r}}=0,
\end{eqnarray}
therefore A can be obtained by
\begin{eqnarray}
&&\nonumber\langle f(x),RB^{T}(x)\rangle_{w_{r}}=\langle A^{T}RB(x),RB^{T}(x)\rangle_{w_{r}},\\
\nonumber\\
\nonumber &&A^{T}=\langle f(x),RB^{T}(x)\rangle_{w_{r}}~\langle RB(x),RB^{T}(x)\rangle_{w_{r}}^{-1},~n=0, 1,\cdots, N.
\end{eqnarray}

\subsection{Exponentioal Bessel Functions}
Exclusive of rational functions we can use exponential transformation to have new functions which are also defined on the semi-infinite interval. The exponential Bessel functions  ($EBs$) can be defined by 
\begin{eqnarray}\label{eqqq06}
&&EB_{n}(x)=B_{n}(1-e^{-x/L})~,~~n=0, 1,\cdots,N.\nonumber
\end{eqnarray}
or
\begin{eqnarray}
&&EB_{n}(x,L)=\sum_{r=0}^{[\frac{N-n}{2}]}\frac{(-1)^r}{r!(n+r)!}(1-e^{-x/L})^{2r+n},~~n=0, 1,\cdots,N.
\end{eqnarray}
where parameter L is a constant parameter and, like rational functions, it sets the length scale of the mapping.

All of the above relations can also be used to EBs with respect to the weight function $w_{e}(x,L)=\frac{e^{-x/L}}{L}$ in the interval $[0,\infty)$.

\section{The quasilinearization method (QLM)}\label{sec3}
The QLM is a generalization of the Newton-Raphson method \cite{QLM02} to solve the nonlinear differential equation as a limit of approximating the nonlinear terms by an iterative sequence of linear expressions. The QLM techniques are based on the linearization of the higher order ordinary/partial differential equation and require the solution of a linear ordinary differential equation at each iteration. Mandelzweig and Tabakin \cite{QLM05} have determined general conditions for the quadratic, monotonic and uniform convergence of the QLM method to solve both initial and boundary value problems in nonlinear ordinary $n$th order differential equations in $N$-dimensional space.

Let us assume that a second-order nonlinear ordinary differential equation in one variable on the interval $[0, \infty)$ as follows:
\begin{eqnarray}\label{lbl12}
\frac{d^2u}{dx^2}=F(u'(x),u(x),x),
\end{eqnarray}
with the boundary conditions: $u(0)=A,~u(\infty)=B$, where $A$ and $B$ are real constants and $F$ is nonlinear function.

By utilizing the QLM to solve Eq. (\ref{lbl12}) determines the $(I+1)$th iterative approximation $u_{I+1}(t)$ as a solution of the linear differential equation:
\begin{equation}\label{lbl13}
\frac{d^2u_{I+1}}{dx^2}=F(u'_r,u_r,x)+(u_{I+1}-u_r)F_u(u'_r,u_r,x)+(u'_{I+1}-u'_r)F_{u'}(u'_r,u_r,x),
\end{equation}
with the boundary conditions:
\begin{equation}\label{lbl14}
u_{I+1}(0)=A,~~~~~u_{I+1}(\infty)=B,
\end{equation}
where $~r=0,1,2,\cdots$ and the functions $F_u = \partial F/ \partial u$ and $F_{u'} = \partial F / \partial u'$ are functional derivatives of functional $F(u'_r,u_r,x)$.

\section{Application of the Methods}\label{sec4}

In this paper, two methods based on RBs collocation method and EBs collocation method for solving Eq. (\ref{lbl4}), with initial conditions of Eq. (\ref{lbl6}), have been considered.

First, by utilizing QLM technique on Eq. (\ref{lbl4}), we have
\begin{eqnarray}\label{lbl15}
xy''_{I+1}(x)+2y'_{I+1}(x)-(m-1)xy^{m}_{I}(x)+mxy_{I+1}(x)y^{m-1}_{I}(x)=0
\end{eqnarray}
with the initial conditions:
\begin{eqnarray}\label{lbl16} 
y_{I+1}(0)=1,~~~~~~y'_{I+1}(0)=0,
\end{eqnarray}
where $I = 0, 1, 2,\cdots .$\\

For rapid convergence is actually enough that the initial guess is sufficiently good to ensure the smallness of just one of the quantity $q_{r} = k||y_{I+1} - y_{I}||$, where $k$ is a constant independent of $I$. Usually, it is advantageous that $y_{0}(t)$
would satisfy at least one of the initial conditions Eq.  (\ref{lbl6})  \cite{QLM10}, thus set
 $y_{0}(x)=1$ for the initial guess of Lane-Emden equation.

Then, we can approximate $y_{I+1}(x)$  by $N+1$ basis of RBs and EBs as follows:
\begin{enumerate}
\item approximating $y_{I+1}(x)$  by $N+1$ basis of RBs:
\begin{enumerate}
\item[] 
\begin{equation}\label{lbl17}
y_{I+1}(x) \thickapprox u_{N,I+1}(x)=1+x^2\sum^{N}_{n=0}{\hat{b_{n}}RB_{n}(x,L)}.
\end{equation}
where $r = 0, 1, 2,\cdots .$ and two terms $1$ and $x^2$ are to satisfy initial conditions Eq. (\ref{lbl16}).

To apply the collocation method, we have constructed the residual function for  $(I+1)$th iteration in QLM by substituting $y_{I+1}(x)$ by $u_{N,I+1}(x)$ into Eq. (\ref{lbl15}) as following:
\small
\begin{eqnarray}\label{lbl18}
\nonumber RESr_{I+1}(x)=xu''_{N+1,I+1}(x)+2u'_{N+1,I+1}(x)\\-(m-1)xu^{m}_{N+1,I}(x)+  mxu_{N+1,I+1}(x)u^{m-1}_{N+1,I}(x)=0
\end{eqnarray}

\end{enumerate}
\item approximating $y_{I+1}(x)$  by $N+1$ basis of EBs:
\item[] 
\begin{equation}\label{lbl19}
y_{I+1}(x) \thickapprox w_{N,I+1}(x)=\frac{1}{x^2+1}+\frac{x^2}{x+1}\sum^{N}_{n=0}{\hat{c_{n}}EB_{n}(x,L)}.
\end{equation}
where $r = 0, 1, 2,\cdots .$ \\
In this paper, we have been considered two terms $\frac{1}{x^2+1}$ and $\frac{x^2}{x+1}$ to satisfy initital conditions Eq. (\ref{lbl16}).\\
Also, like above, to apply the collocation method, we have constructed the residual function for  $(I+1)$th iteration in QLM by substituting $y_{I+1}(x)$ by $w_{N,I+1}(x)$ into Eq. (\ref{lbl15}) as following:
\begin{eqnarray}\label{lbl20}
\nonumber RESe_{I+1}(x)=xw''_{N+1,I+1}(x)+2w'_{N+1,I+1}(x)\\-(m-1)xw^{m}_{N+1,I}(x)+  mxw_{N+1,I+1}(x)w^{m-1}_{N+1,I}(x)=0
\end{eqnarray}
\end{enumerate}
In all of the spectral methods, the purpose is to find $\hat{b_{n}}$ and $\hat{c_{n}}$  coefficients.

A method for forcing the residual functions Eq. (\ref{lbl18}) and Eq. (\ref{lbl20}) to zero can be defined as collocation algorithm. There is no limitation to choose the point in collocation method. The $N+1$ collocation points have been substituted in $RESr_{I+1}(x)$ and $RESe_{I+1}(x)$ equations, therefore:
\begin{eqnarray}\label{eqqq18}
RESr_{I+1}(x_{i})=0,~~i=0, 1, ,\cdots, N+1.\\
RESe_{I+1}(x_{i})=0,~~i=0, 1, ,\cdots, N+1.
\end{eqnarray}
which $x_{i}$ are roots of the shifted Chebyshev functions on finite interval\cite{refdd27}. Finally, a linear system of equations has been obtained, all of these equations can be solved by Newton method for the unknown coefficients.

\section{Results and discussion}\label{sec5}
The Lane-Emden type describe the variation of density as a function of the radial distance for a polytrope. They was first studied by the astrophysicists Jonathan Homer Lane and Robert Emden, which considered the thermal behavior of a spherical cloud of gas acting under the mutual attraction of its molecules and subject to the classical laws of thermodynamics\cite{refff01,refff02}.  In the Lane-Emden type equations, the first zero of $y(x)$ is an important point of the function, so we have computed $y(x)$ to this zero. In this paper, the equation is solved for $m=1.5, 2, 2.5 ,3$ and $4$, which does not have exact solutions.\\
The comparison of the initial slope $y'(0)$ calculated by RBs-QLM with values obtained by Horedt\cite{refff05} is given in table 1.\\
Table 2 and 3 have presented some numerical examples to illustrate the accuracy and convergence of our suggested methods by increasing the number of points and iterations.\\
Tables 4, 5, 6, 7 and 8 show the obtained values of $y(x)$ and $y'(x)$ by the approch which based on RBs collocation method, for $m=1.5, 2, 2.5 ,3$ and $4$ with values of $N=75$ and iteration 15.\\
The resulting graphs of the standard Lane-Emden equation  obtained by present methods for  $m=1.5, 2, 2.5 ,3$ and $4$ are shown in figure 1.\\
 Finally, figures 2, 3, 4, 5 and 6 show the residual errors for approximation solutions by basis of rational and exponential function with $N=50, 75$ and 100. note that the residual error decreases with the increase of the collocation points.
\section{Conclusion}\label{sec6}
The fundamental goal of this paper were to introduce novel hybrid basis of Rational Bessel and Exponential functions (RB, EBs) with quasilinearization meyhod (QLM) to construct an approximation for solving nonlinear Lane-Emden type equations. These problems describe a variety of phenomena in theoretical physics and astrophysics, including aspects of stellar structure, the thermal history of a spherical cloud of gas, isothermal gas spheres, and thermionic currents \cite{refff04}. To achieve these goal at first, a sequence of linear differential equations is obtained by utilizing QLM over Lane-Emden equation. Second, in each iteration of QLM, the linear differential equation is solved by new RBs and EBs collocation method. This paper has been shown that the present works have provided two acceptable approaches for solving Lane-Emden type equations coused by the following reasons: 
\begin{itemize}
  \item[1]  Cause of simplicity to solve problems and convergence of approximation functions, we convert the nonlinear problems to a sequence linear equations using QLM
  \item[2]  Numerical results indicate effectiveness, applicability and accuracy of the present approaches.
  \item[3] Present paper described shortly bibliography of different methods utilizined previous works
for solving Lane-Emden-type equations.
  \item[4] The approaches applied to solve the problems without reformulating the equation to bounded domains.
  \item[5] The approaches have been displyed converges when increasing the number of collocation points by tabular reports.
\item[6] At the first time, Rational Bessel and Exponential functions have been to obtain numerical outcomes of the nonlinear exponent $m$ of the standard Lane-Emden equations.
\item[7] Moreover, a very good approximation solution of $y(x)$ for Lane-Emden type equations with various value of parameter $m$ after only fifteen iterations are obtained. So, these methods are a good experience and method for the other sciences.

\end{itemize}

.\\\\\\\\\\\\\\\\\\\\\\\\\\\\\\\\\\\\\\\\\\
%%%%%%%%%%%%%%%%%%%%%%%%%%%%%%%%%%%%%%%%%%
\begin{figure}[!ht]
\centering
\subfigure[rational]{
\includegraphics*[width=5.5cm]{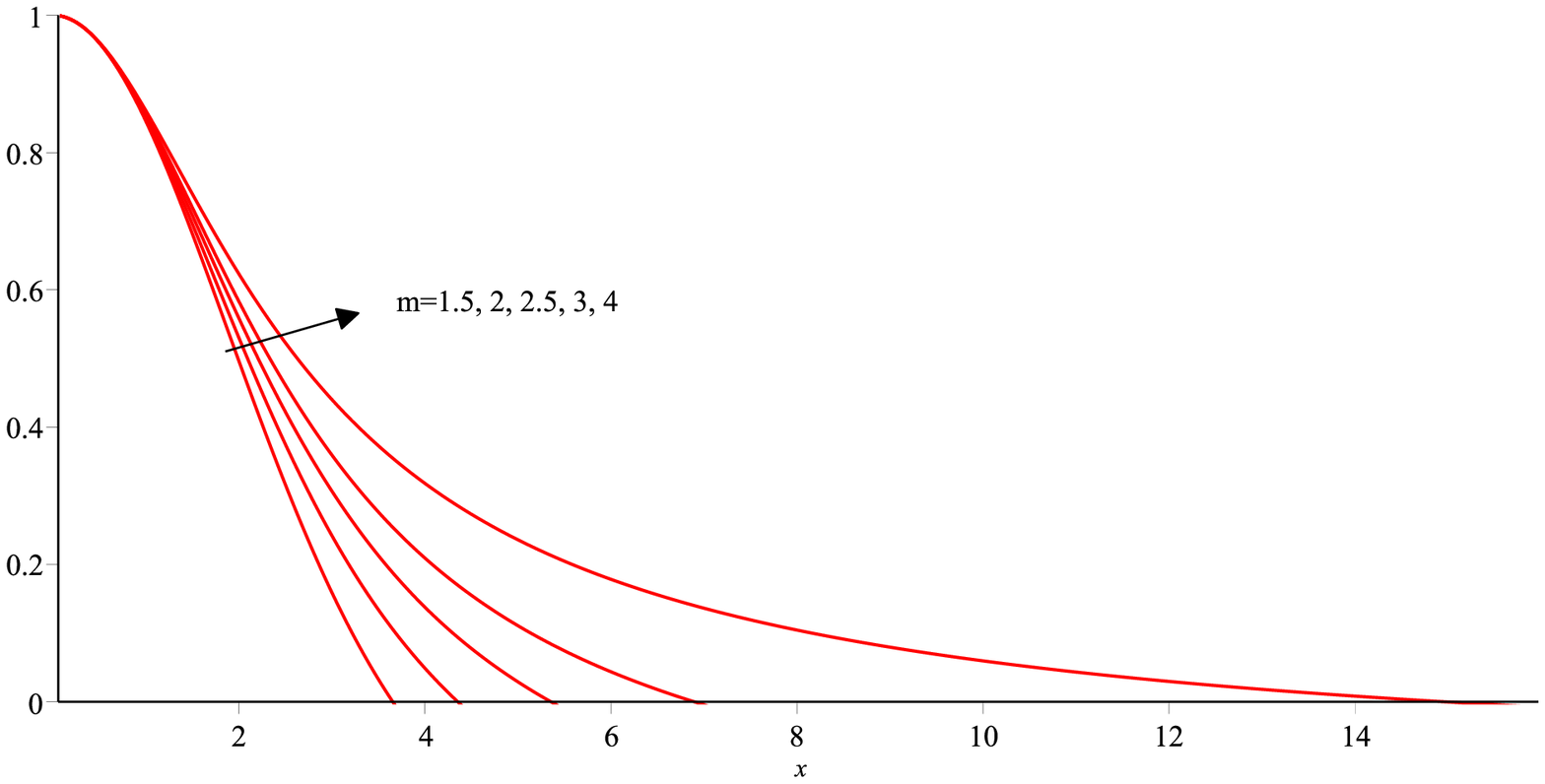}}
\hspace{3mm}
\subfigure[exponential]{
\includegraphics*[width=5.5cm]{Y.eps}}
\caption{The obtained graphs of solutions of Lane-Emden standard equations by basis of  RBs and EBs with $m=1.5, 2, 2.5, 3, 4$ }

\end{figure}

%%%%%%%%%%%%%%%%%%%%%%%%%%%%%%%%%%%%%%%%%%

%%%%%%%%%%%%%%%%%%%%%%%%%%%%%%%%%%%%%%%%%%

\begin{figure}[!ht]
\centering
\subfigure[rational]{
\includegraphics*[width=5cm]{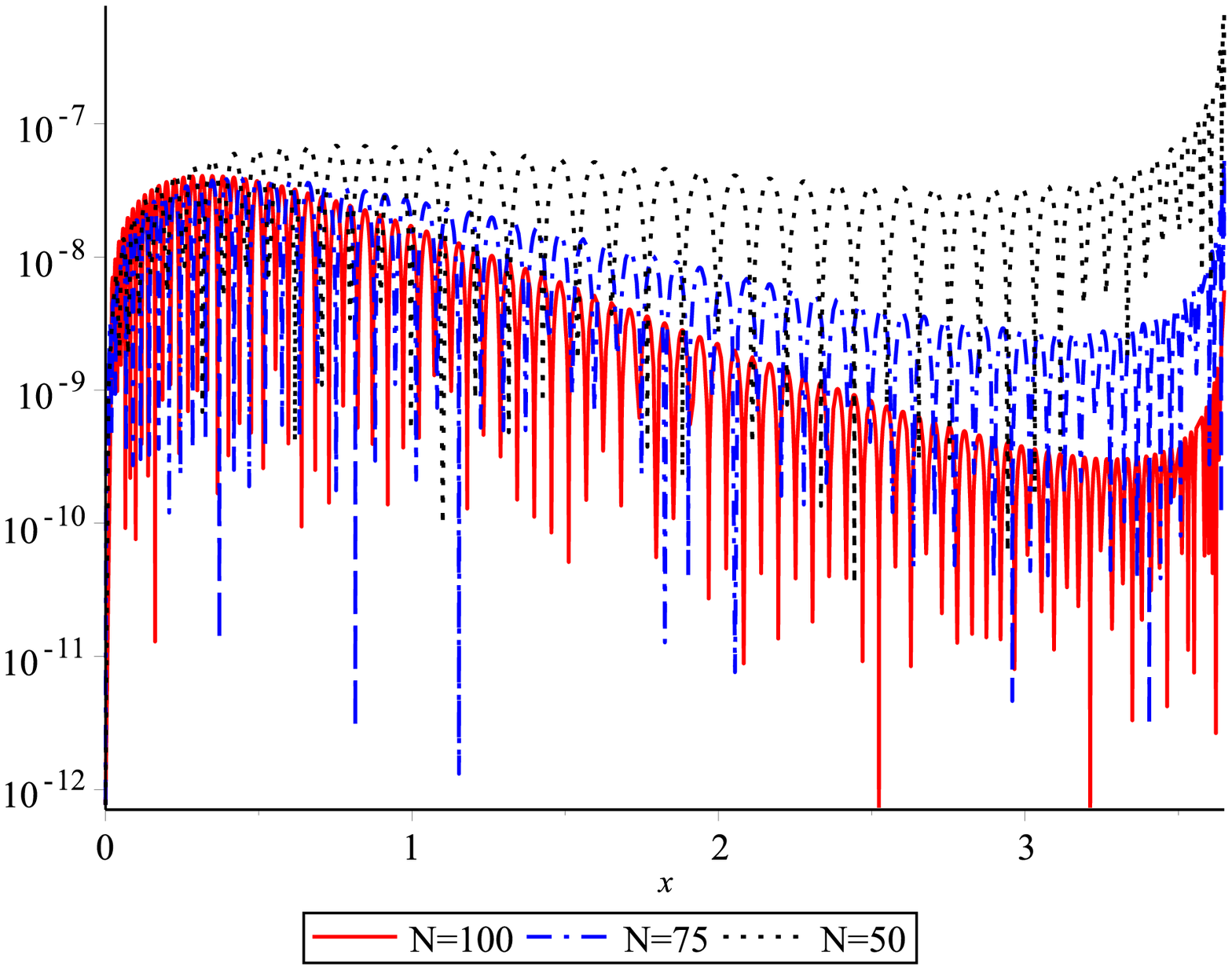}}
\hspace{3mm}
\subfigure[exponential]{
\includegraphics*[width=5cm]{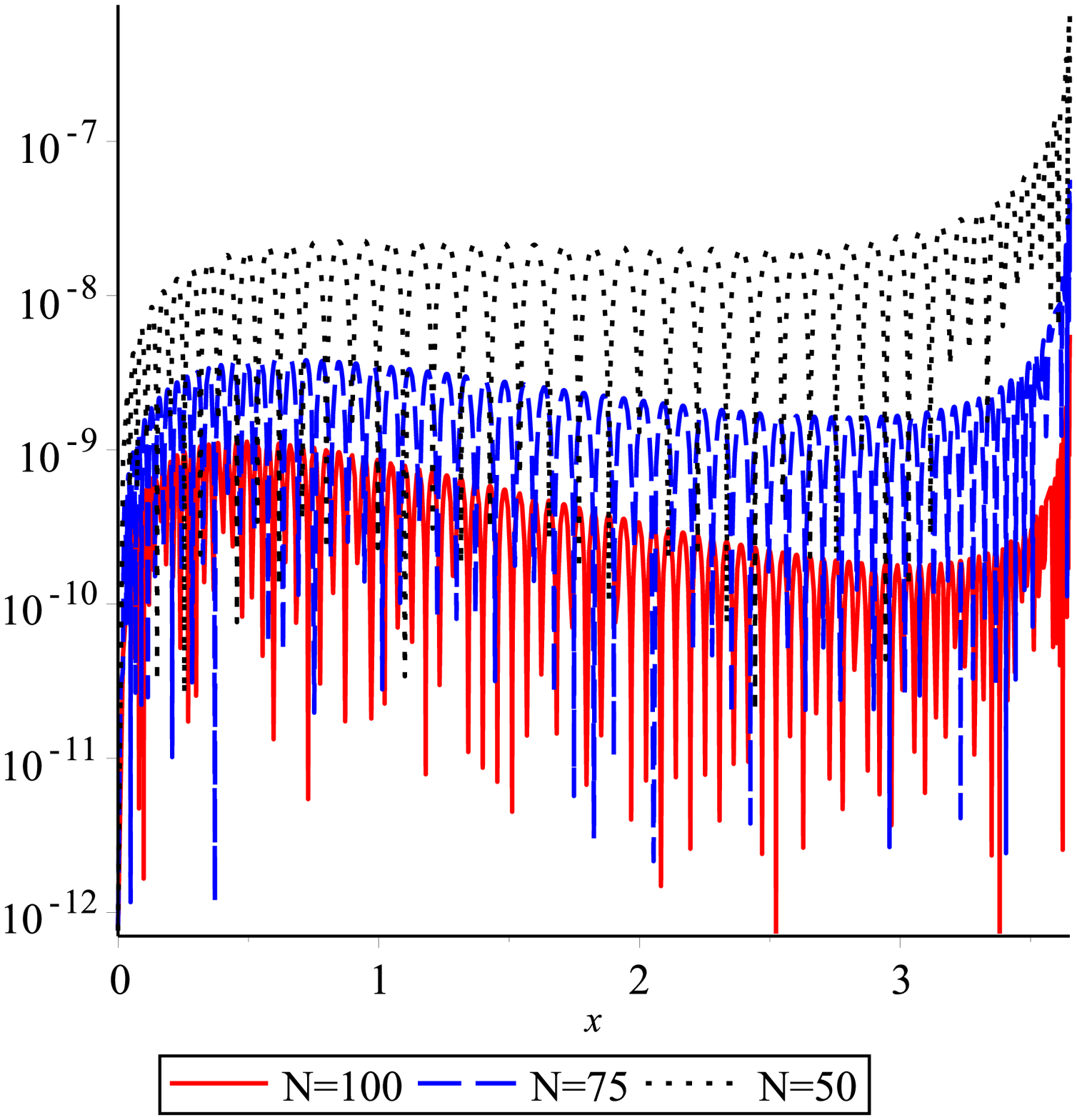}}
\caption{Logarithmic graph of residual error by present works with $N=50,75,100$ and iteration 15 when m=1.5.}

\end{figure}

%%%%%%%%%%%%%%%%%%%%%%%%%%%%%%%%%%%%%%%%%%

%%%%%%%%%%%%%%%%%%%%%%%%%%%%%%%%%%%%%%%%%%

\begin{figure}[!ht]
\centering
\subfigure[rational]{
\includegraphics*[width=5cm]{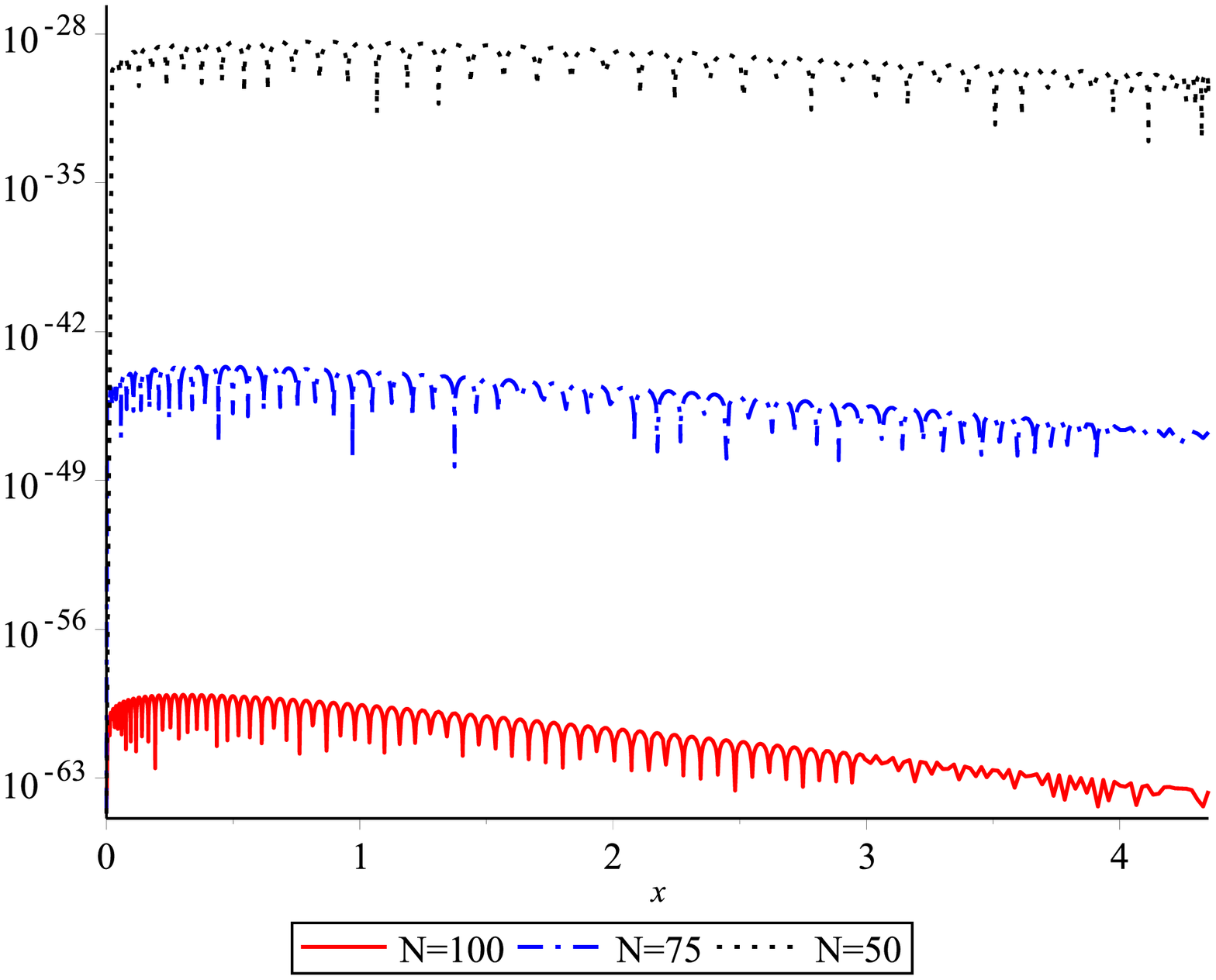}}
\hspace{3mm}
\subfigure[exponential]{
\includegraphics*[width=5cm]{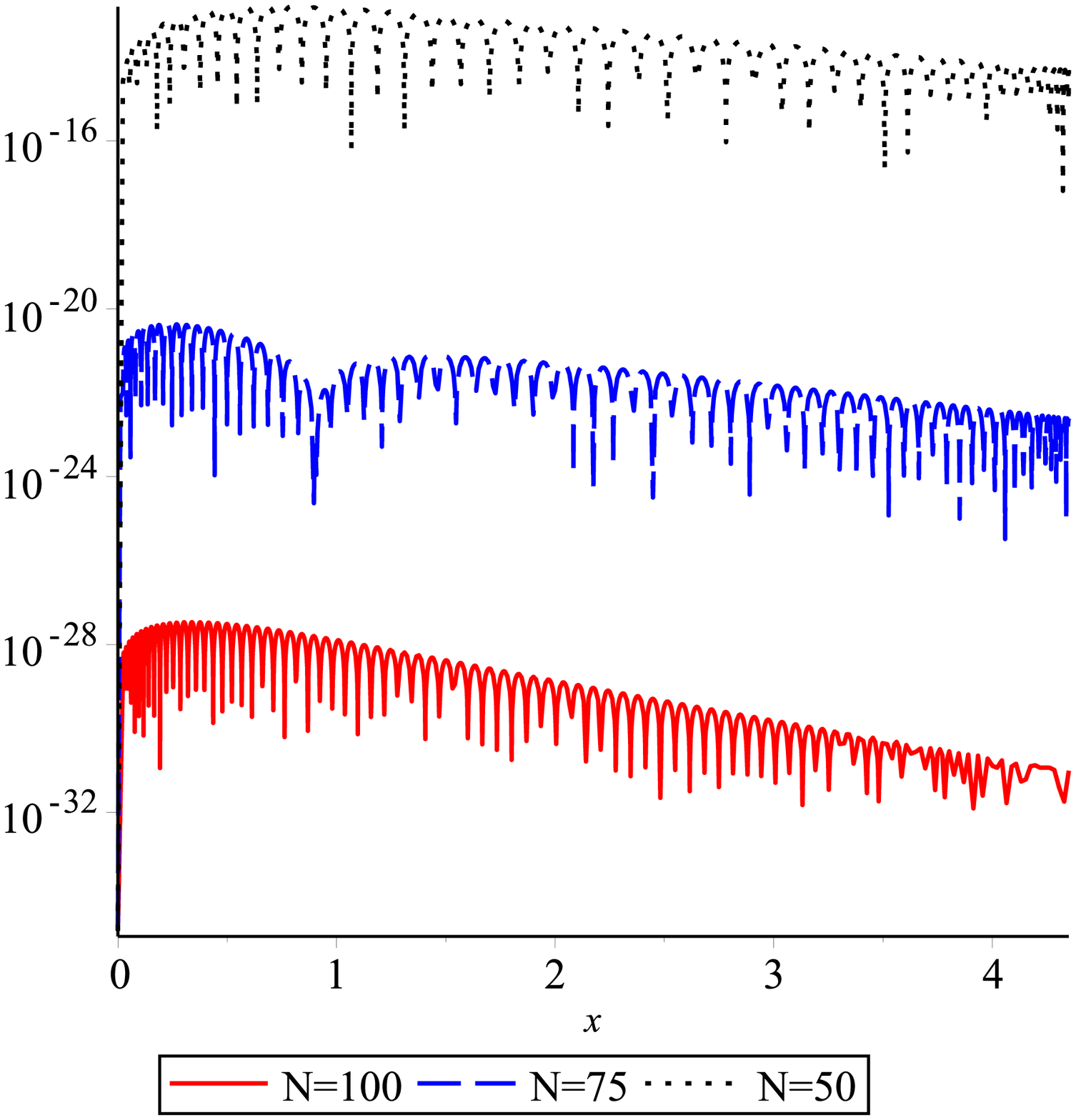}}
\caption{Logarithmic graph of residual error by present works with $N=50,75,100$ and iteration 15 when m=2.}
\end{figure}

%%%%%%%%%%%%%%%%%%%%%%%%%%%%%%%%%%%%%%%%%%

%%%%%%%%%%%%%%%%%%%%%%%%%%%%%%%%%%%%%%%%%%

\begin{figure}[!ht]
\centering
\subfigure[rational]{
\includegraphics*[width=5cm]{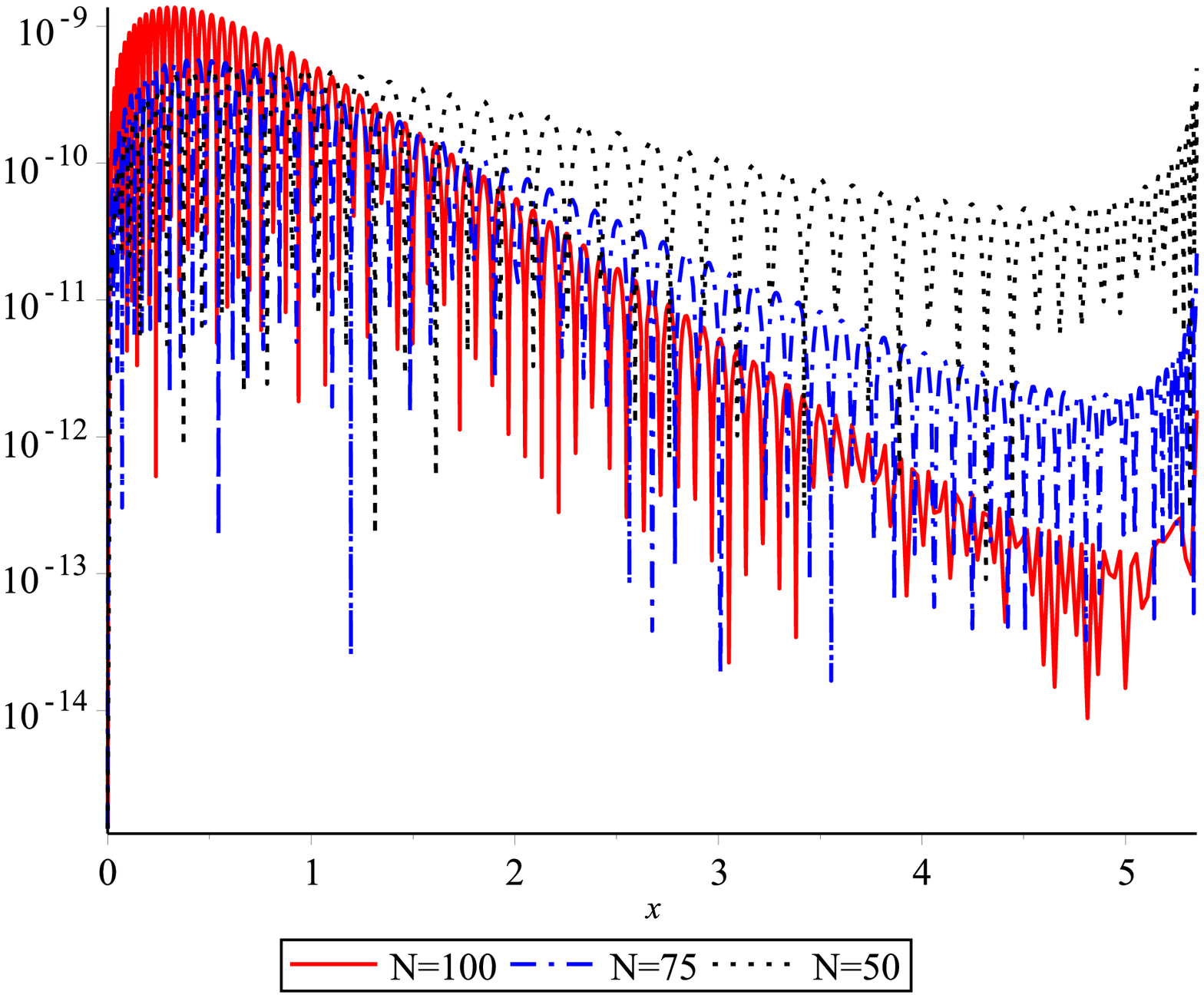}}
\hspace{3mm}
\subfigure[exponential]{
\includegraphics*[width=5cm]{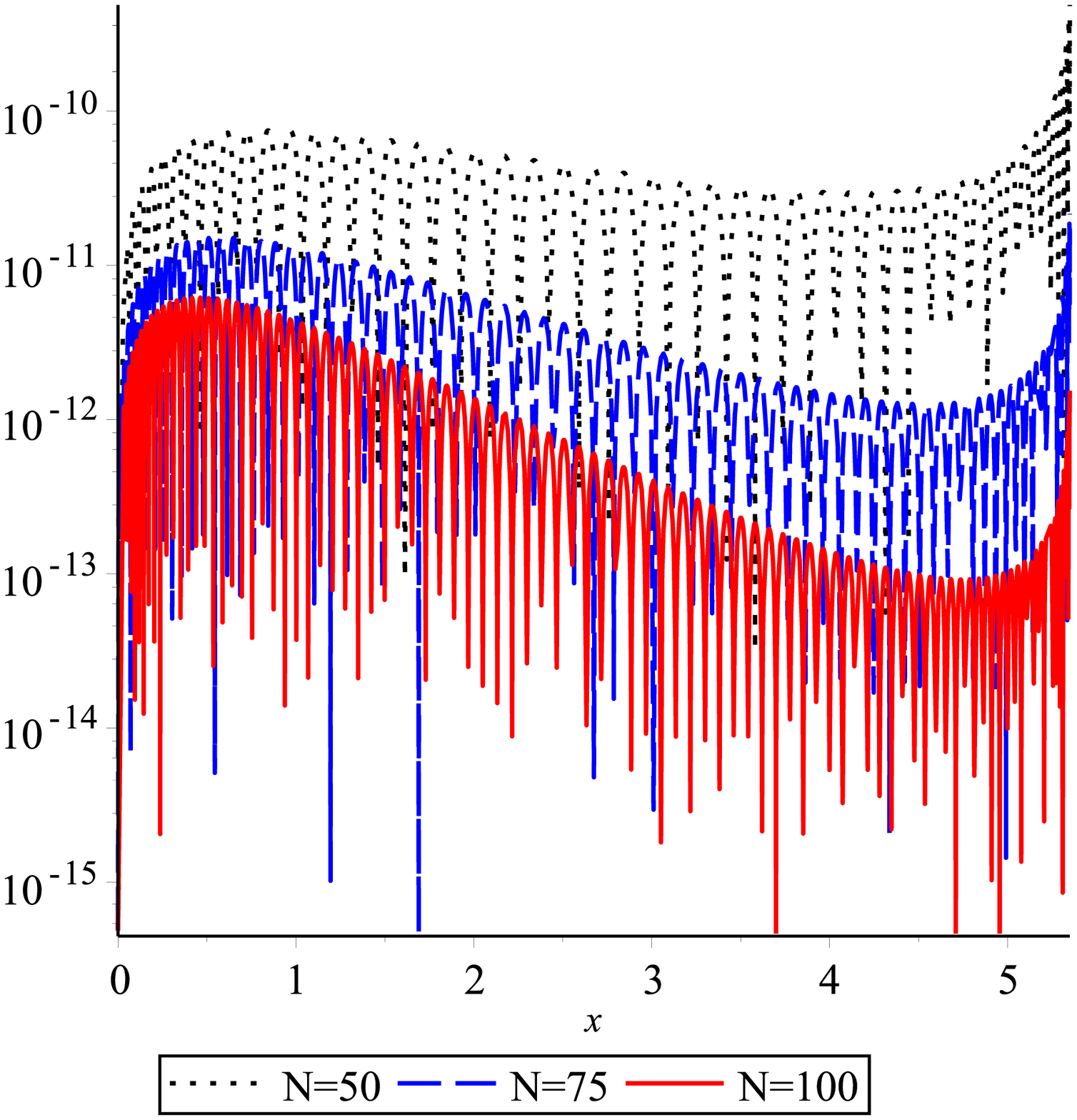}}
\caption{Logarithmic graph of residual error by present works with $N=50,75,100$ and iteration 15 when m=2.5.}
\end{figure}

%%%%%%%%%%%%%%%%%%%%%%%%%%%%%%%%%%%%%%%%%%

%%%%%%%%%%%%%%%%%%%%%%%%%%%%%%%%%%%%%%%%%%

\begin{figure}[!ht]
\centering
\subfigure[rational]{
\includegraphics*[width=5cm]{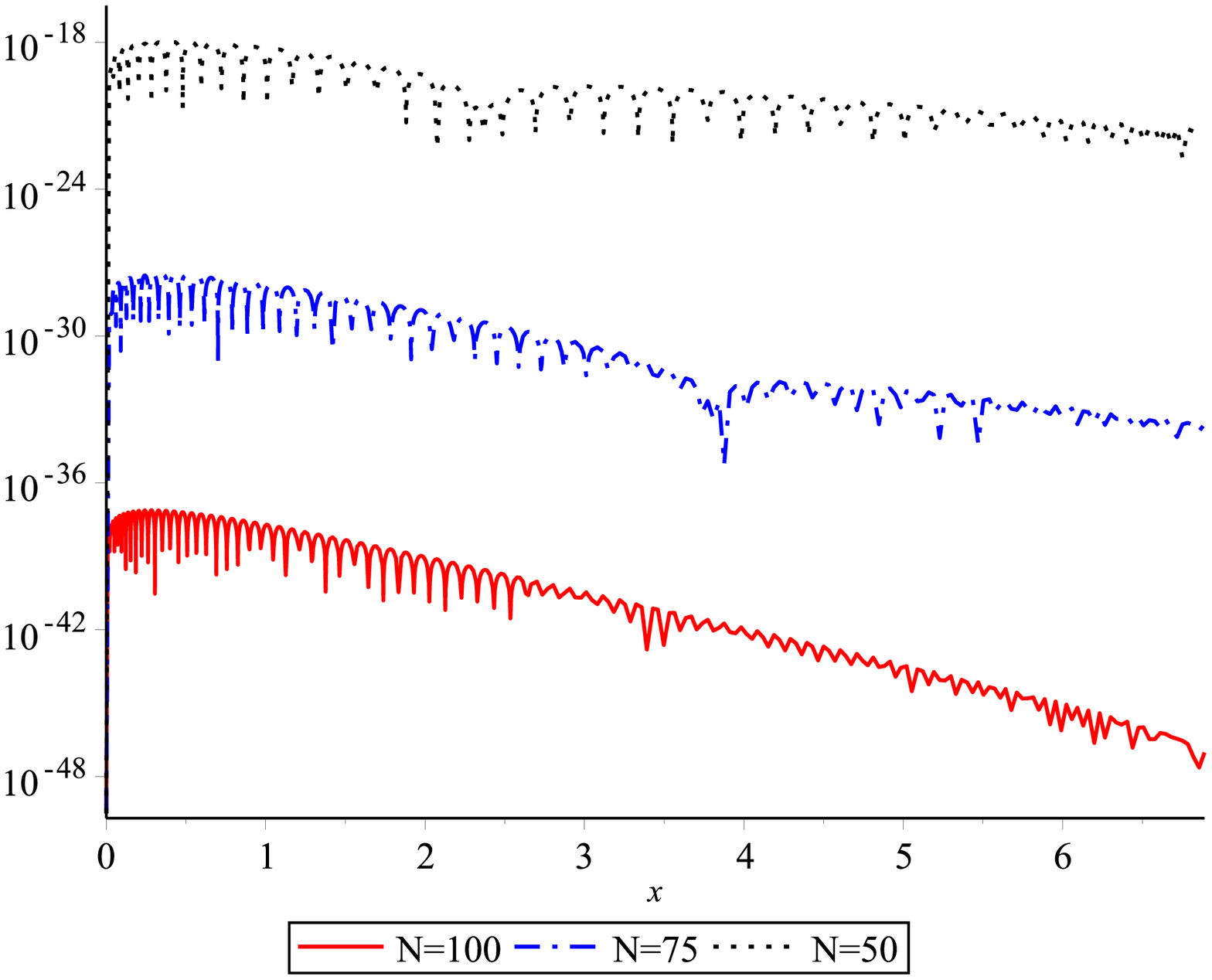}}
\hspace{3mm}
\subfigure[exponential]{
\includegraphics*[width=5cm]{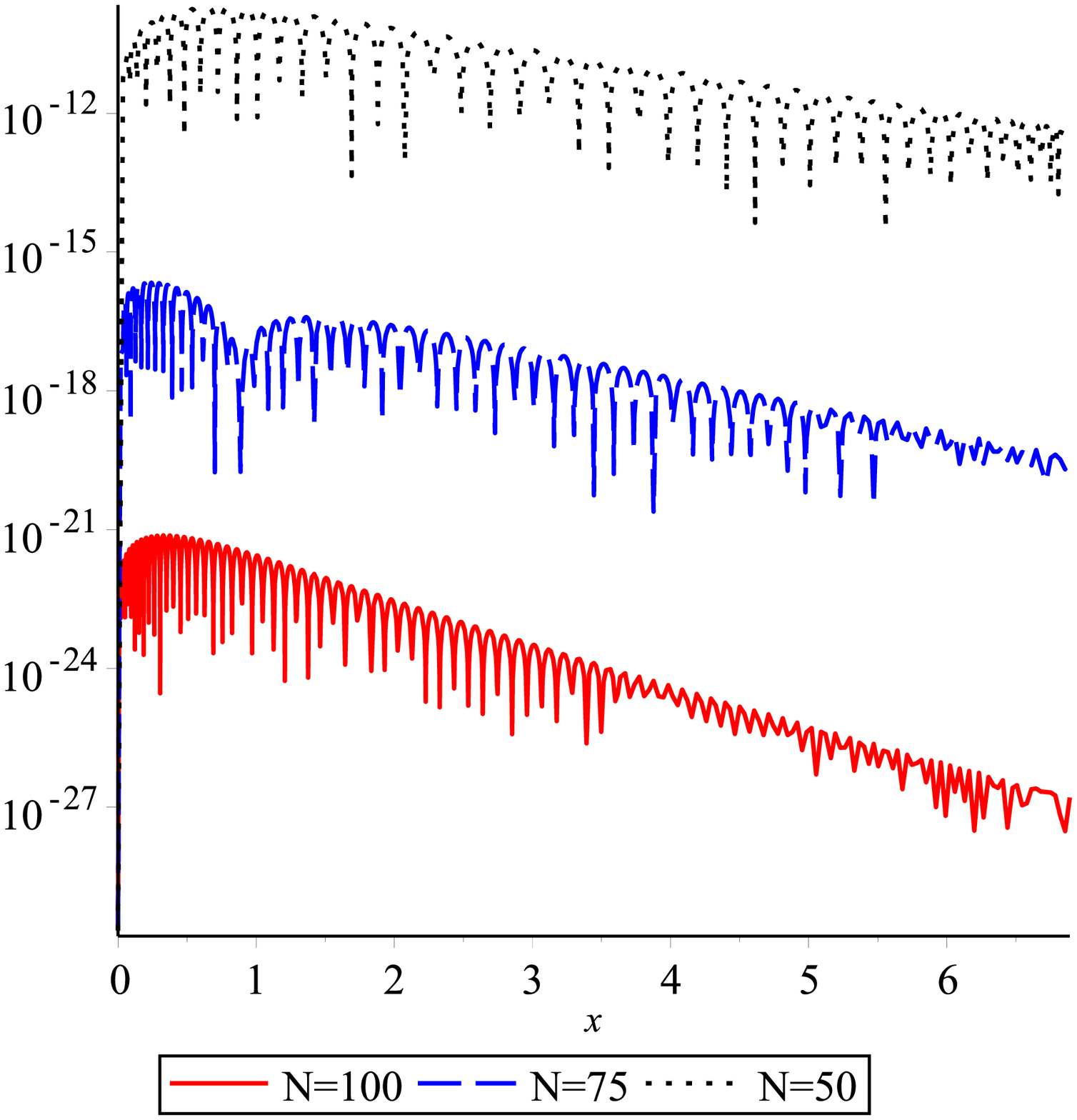}}
\caption{Logarithmic graph of residual error by present works with $N=50,75,100$ and iteration 15 when m=3.}
\end{figure}

%%%%%%%%%%%%%%%%%%%%%%%%%%%%%%%%%%%%%%%%%%

%%%%%%%%%%%%%%%%%%%%%%%%%%%%%%%%%%%%%%%%%%

\begin{figure}[!ht]
\centering
\subfigure[rational]{
\includegraphics*[width=5cm]{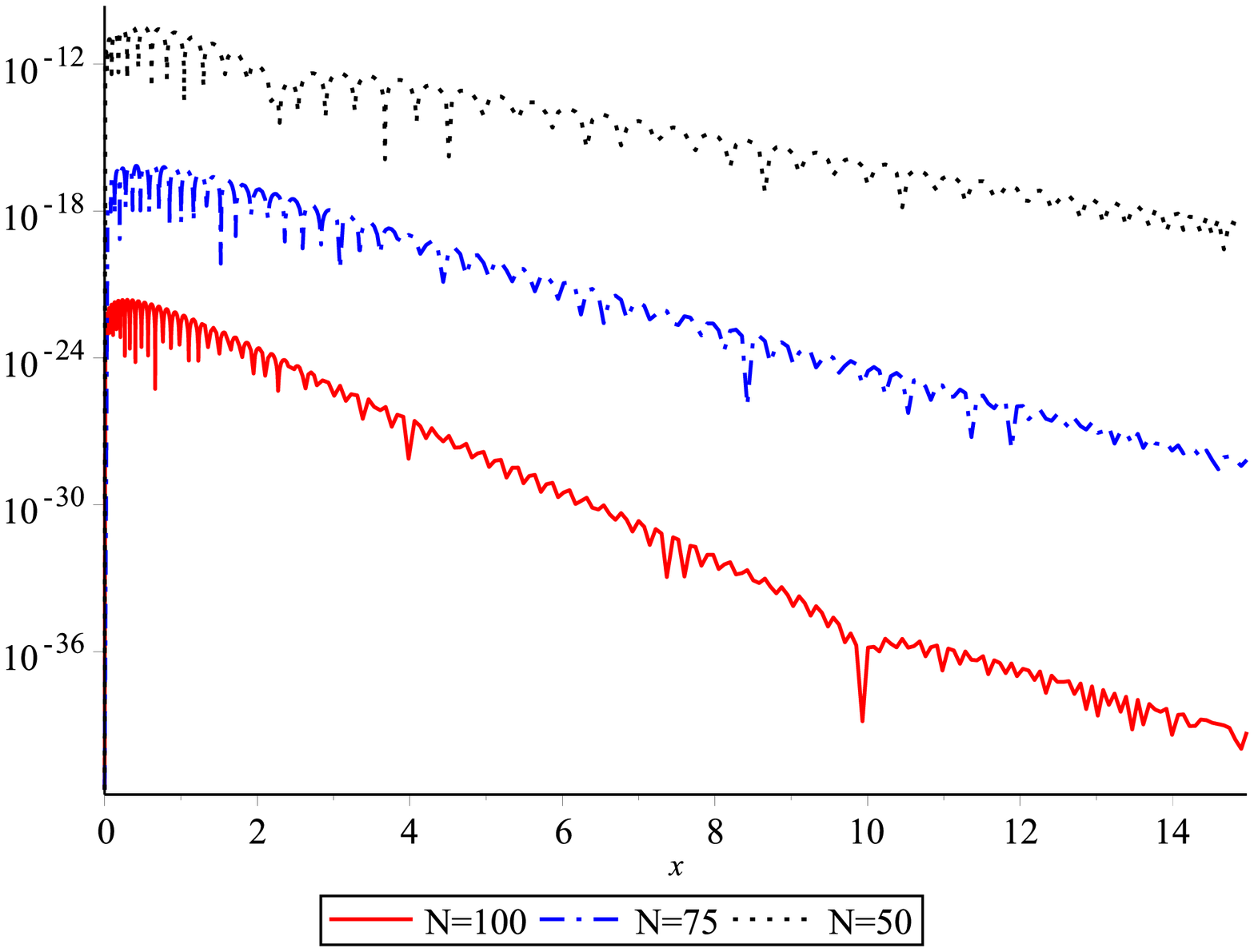}}
\hspace{3mm}
\subfigure[exponential]{
\includegraphics*[width=5cm]{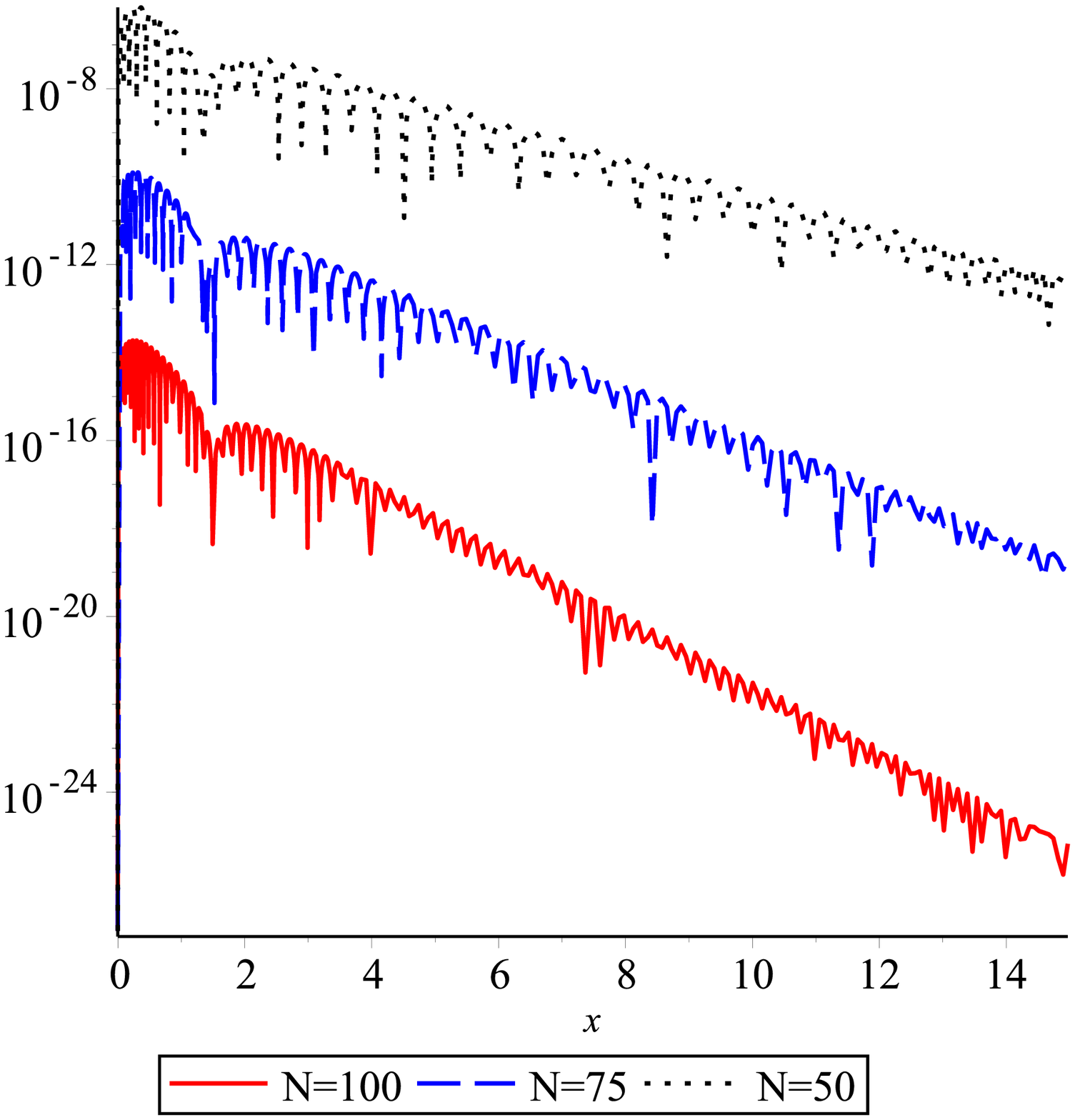}}
\caption{Logarithmic graph of residual error by present works with $N=50,75,100$ and iteration 15 when m=4.}
\end{figure}

%%%%%%%%%%%%%%%%%%%%%%%%%%%%%%%%%%%%%%%%%%

%%%%%%%%%%%%%%%%%%%%%%%%%%%%%%%%%%%%%%%%%%
\tiny  
\begin{center}
\begin{table}[ht!]
\caption{Comparison of the first zeros of standard Lane-Emden equations, for numerical values given by Horedt\cite{refff05}  nad the present methods with  N=75 and iteration 15}
\centering
%\footnotesize 
\begin{tabular}[ht!]{c|l|l|l}
\toprule
   $m$  & RBs & EBs& Horedt \cite{refff05}    \\ 
\midrule       
     1.5   & 3.65375373622763424836747856706295570  &  3.653753736227530116708951&3.65375374     \\ 
    2.0   &  4.35287459594612467697357006152614262  &  4.352874595946124676973570&4.35287460   \\ 
    2.5    & 5.35527545901076012377857991160851840  &  5.355275459010769844745925&5.35527546  \\ 
    3.0    & 6.89684861937696037545452818712314053  &  6.896848619376960375436984&6.89684862   \\ 
    4.0    & 14.9715463488380950976509645543077611  & 14.97154634883796085494984&14.9715463  \\ 

\bottomrule
\end{tabular}
\end{table}
\end{center}
%%%%%%%%%%%%%%%%%%%%%%%%%%%%%%%%%%%%%%%%%%

%%%%%%%%%%%%%%%%%%%%%%%%%%%%%%%%%%%%%%%%%%
\footnotesize  
\begin{center}
\begin{table}[ht!]
\caption{Numerical results of first zeros by basis of RBs with various values of $m$, $N$ and iterations, accurate digits are bold }
\centering
%\footnotesize 
\begin{tabular}[ht!]{cll|l}
\toprule
   $m$ &$N$&  iteration  & RBs    \\

%%%%%%%%%%%%%%%%%%%%%%%%%%%%%%%%%%%%%%%%%%%%%%%%% m=1.5
\midrule       
     1.5& 50 &05  &  \textbf{3.6537537362}5072342590     \\ 
          &     &10  &  \textbf{3.6537537362}5071853754\\
          &     &15  &  \textbf{3.6537537362}5071853754  \\
          &     &20  &  \textbf{3.6537537362}5071853754\\ 
\\
          &75 &05  &  \textbf{3.65375373622}763914172     \\ 
          &     &10  & \textbf{3.65375373622}763424836\\
          &     &15  & \textbf{3.65375373622}763424836  \\
          &     &20  & \textbf{3.65375373622}763424836 \\
\\
         &100&05  &  \textbf{3.65375373622}225950682     \\ 
          &    &10  &  \textbf{3.65375373622}225461061\\
          &    &15  &  \textbf{3.65375373622}225461061  \\
          &    &20  &  \textbf{3.65375373622}225461061 \\

%%%%%%%%%%%%%%%%%%%%%%%%%%%%%%%%%%%%%%%%%%%%%%%%% m=2
\midrule       
        2& 50&05  &   \textbf{4.35}41023191782544510394699271974639349588062470049419121696397470     \\ 
          &    &10  &   \textbf{4.35287459594612467697357006152614}339487342457587311708331752\\
          &    &15  &   \textbf{4.35287459594612467697357006152614}339487342457587311708331752  \\
          &    &20  &   \textbf{4.35287459594612467697357006152614}339487342457587311708331752\\ 
\\
          &75 &05  &  \textbf{4.35287459}7893199784546816142774753394907169932534281348066892095     \\ 
          &     &10  &  \textbf{4.352874595946124676973570061526142628112365363213}147181521\\
          &     &15  &  \textbf{4.352874595946124676973570061526142628112365363213}147181521 \\
          &     &20  &  \textbf{4.352874595946124676973570061526142628112365363213}147181521 \\
\\
          &100&05  &  \textbf{4.35287459}7893199784546816142774753394907169932542963806389373524     \\ 
          &     &10  &   \textbf{4.352874595946124676973570061526142628112365363213}008835302\\
          &     &15  &   \textbf{4.352874595946124676973570061526142628112365363213}008835302 \\
          &     &20  &   \textbf{4.352874595946124676973570061526142628112365363213}008835302\\
%%%%%%%%%%%%%%%%%%%%%%%%%%%%%%%%%%%%%%%%%%%%%%%%% m=2.5
\midrule       
     2.5& 50&05  &  \textbf{5.35529645}45076443677     \\ 
          &    &10  &  \textbf{5.3552754590107}44925\\
          &    &15  &  \textbf{5.3552754590107}44925  \\
          &    &20  &  \textbf{5.3552754590107}44925 \\ 
\\
         &75 &05  &  \textbf{5.3552}9645450764436772      \\ 
          &   &10  &   \textbf{5.3552754590107}601237\\
          &   &15  &   \textbf{5.3552754590107}601237  \\
          &   &20  &   \textbf{5.3552754590107}601237 \\
\\
         &100&05&  \textbf{5.3552}9645450764436772     \\ 
         &     &10&   \textbf{5.3552754590107}873176\\
         &     &15&  \textbf{5.3552754590107}873176  \\
         &     &20&  \textbf{5.3552754590107}873176 \\

%%%%%%%%%%%%%%%%%%%%%%%%%%%%%%%%%%%%%%%%%%%%%%%%% m=3
\midrule       
        3& 50 &05  &  7.1216938046517305045330727094680858444666907392     \\ 
          &     &10  &  \textbf{6.896848619376960375454}2796110144170369244612\\
          &     &15  &  \textbf{6.896848619376960375454}2796110144170369244612  \\
          &     &20  &  \textbf{6.896848619376960375454}2796110144170369244612 \\ 
\\
          & 75&05  &  7.1216938046404145204995503800811081360235860196     \\ 
          &     &10  &  \textbf{6.8968486193769603754545281871231}4053555203\\
          &     &15  &  \textbf{6.8968486193769603754545281871231}4053555203  \\
          &     &20  &  \textbf{6.8968486193769603754545281871231}4053555203 \\
\\
          &100&05  &  7.1216938046404152911963760032858519494248670403     \\ 
          &     &10  &   \textbf{6.8968486193769603754545281871231}2127697218\\
          &     &15  &   \textbf{6.8968486193769603754545281871231}2127697218  \\
          &     &20  &   \textbf{6.8968486193769603754545281871231}2127697218 \\

%%%%%%%%%%%%%%%%%%%%%%%%%%%%%%%%%%%%%%%%%%%%%%%%% m=4
\midrule       
        4& 50 &05  &  \textbf{1}6.711045707072842315340457798698905988740559701     \\ 
          &     &10  &   \textbf{14.97154}867059731700938111496437106672775015032\\
          &     &15  &   \textbf{14.971546348838}2089901971898981867391578167481  \\
          &     &20  &   \textbf{14.971546348838}2089901971898981867391578167481 \\ 
\\
          &75 &05  &   16.402670239960775259418702056564527058250944781     \\ 
          &    &10  &    \textbf{14.97154}289318059650158197244640609252173187180\\
          &    &15  &    \textbf{14.9715463488380950976}50964554307761107155441 \\
          &    &20  &    \textbf{14.9715463488380950976}50964554307761107155441 \\
\\
         &100&05  &   16.172787459355139190211994543646969560813181439     \\ 
         &     &10  &    \textbf{14.97154}439717955256111887952830248179390503419\\
         &     &15  &    \textbf{14.9715463488380950976}11066133148254587457821 \\
         &     &20  &    \textbf{14.9715463488380950976}11066133148254587457821\\

\bottomrule
\end{tabular}
\end{table}
\end{center}
%%%%%%%%%%%%%%%%%%%%%%%%%%%%%%%%%%%%%%%%%%

%%%%%%%%%%%%%%%%%%%%%%%%%%%%%%%%%%%%%%%%%%
\footnotesize 
\begin{center}
\begin{table}[ht!]
\caption{Numerical results of first zeros by basis of EBs with various values of $m$, $N$ and iterations, accurate digits are bold  }
\centering
%\footnotesize 
\begin{tabular}[ht!]{cll|l}
\toprule
   $m$ &$N$&  iteration  & EBs    \\

%%%%%%%%%%%%%%%%%%%%%%%%%%%%%%%%%%%%%%%%%%%%%%%%% m=1.5
\midrule       
      1.5& 50 &05  &  \textbf{3.6537537362}5083916424     \\ 
          &      &10  &   \textbf{3.6537537362}5083427589\\
          &      &15  &   \textbf{3.6537537362}5083427589  \\
          &      &20  &   \textbf{3.6537537362}5083427589 \\ 
\\
          &75 &05  &    \textbf{3.65375373622}753501010     \\ 
          &     &10  &    \textbf{3.65375373622}753011670\\
          &     &15  &   \textbf{3.65375373622}753011670  \\
          &     &20  &   \textbf{3.65375373622}753011670 \\
\\
         &100&05  &  \textbf{3.65375373622}227432714     \\ 
          &    &10  &   \textbf{3.65375373622}226943093\\
          &    &15  &   \textbf{3.65375373622}226943093  \\
          &    &20  &   \textbf{3.65375373622}226943093 \\

%%%%%%%%%%%%%%%%%%%%%%%%%%%%%%%%%%%%%%%%%%%%%%%%% m=2
\midrule       
         2& 50 &05 &  \textbf{4.35287459}7893199785338903310594652764    \\ 
          &      &10  &  \textbf{4.35287459594612467}776565735834309221\\
          &      &15  &  \textbf{4.35287459594612467}776565735834309221 \\
          &      &20  &  \textbf{4.35287459594612467}776565735834309221 \\ 
\\
          & 75 &05  & \textbf{4.35287459594612467697}472244822039342    \\ 
          &     &10  &  \textbf{4.352874595946124676973570}1033024306\\
          &     &15  &  \textbf{4.352874595946124676973570}1033024306  \\
          &     &20  &  \textbf{4.352874595946124676973570}1033024306 \\
\\
         &100&05 &  \textbf{4.35287459}78931997845468161427747526020    \\ 
          &    &10  &  \textbf{4.352874595946124676973570}0615261418 \\
          &    &15  &  \textbf{4.352874595946124676973570}0615261418   \\
          &    &20  &  \textbf{4.352874595946124676973570}0615261418 \\
%%%%%%%%%%%%%%%%%%%%%%%%%%%%%%%%%%%%%%%%%%%%%%%%% m=2.5
\midrule       
     2.5& 50 &05  & \textbf{5.35529645}450764438595     \\ 
          &     &10  &  \textbf{5.3552754590107}203902\\
          &     &15  &  \textbf{5.3552754590107}203902  \\
          &     &20  &  \textbf{5.3552754590107}203902 \\ 
\\
          & 75&05  & \textbf{5.35529645}450764436772     \\ 
          &    &10  &  \textbf{5.3552754590107}698447\\
          &    &15  &  \textbf{5.3552754590107}698447  \\
          &    &20  &  \textbf{5.3552754590107}698447 \\
\\
         &100&05  & \textbf{5.35529645}450764436772    \\ 
          &    &10  &  \textbf{5.3552754590107}770840\\
          &    &15  &  \textbf{5.3552754590107}770840  \\
          &    &20  &  \textbf{5.3552754590107}770840 \\

%%%%%%%%%%%%%%%%%%%%%%%%%%%%%%%%%%%%%%%%%%%%%%%%% m=3
\midrule       
        3& 50 &05  &   7.12169371888993111013538427437823   \\ 
          &     &10  &   \textbf{6.89684861937696}9505160794512467\\
          &     &15   &  \textbf{6.89684861937696}9505160794512467  \\
          &     &20   & \textbf{6.89684861937696}9505160794512467\\ 
\\
         & 75 &05  &    7.12169380466912339539903047482119   \\ 
          &     &10  &  \textbf{6.8968486193769603754}3698467213 \\
          &     &15  &  \textbf{6.8968486193769603754}3698467213  \\
          &     &20  &  \textbf{6.8968486193769603754}3698467213 \\
\\
          &100&05  & \textbf{6.89684861937696037}791871227973  \\ 
          &     &10  &  \textbf{6.8968486193769603754}5452817312\\
          &     &15  &  \textbf{6.8968486193769603754}5452817312  \\
          &     &20  &  \textbf{6.8968486193769603754}5452817312 \\

%%%%%%%%%%%%%%%%%%%%%%%%%%%%%%%%%%%%%%%%%%%%%%%%% m=4
\midrule       
        4& 50 &05  &  \textbf{1}6.26491731190237369943385    \\ 
          &     &10  &  \textbf{14.97154}73275763026931076\\
          &     &15  &  \textbf{14.9715463}522353010587855  \\
          &     &20  &  \textbf{14.9715463}522353010587855\\ 
\\
          &75 &05  &  \textbf{1}6.05210011924457026115446     \\ 
          &    &10  &   \textbf{14.97154}72743172097800824\\
          &    &15  &   \textbf{14.97154634883}7960854949  \\
          &    &20  &   \textbf{14.97154634883}7960854949\\
\\
          &100&05  &  \textbf{1}6.03218609785456527010395      \\ 
          &     &10  &  \textbf{14.97154}72744622920651685\\
          &     &15  &  \textbf{14.97154634883}8095104708  \\
          &     &20  &  \textbf{14.97154634883}8095104708\\

\bottomrule
\end{tabular}
\end{table}
\end{center}
%%%%%%%%%%%%%%%%%%%%%%%%%%%%%%%%%%%%%%%%%%

%%%%%%%%%%%%%%%%%%%%%%%%%%%%%%%%%%%%%%%%%%
\footnotesize 
\begin{center}
\begin{table}[ht!]
\caption{Obtained values of $y(x)$ and $y'(x)$ of  standard Lane-Emden equations for m = 1.5 by basis of RBs  with $N=75$ and iterations 15 }
\centering
%\footnotesize 
\begin{tabular}[ht!]{c|l|l}
\toprule
   $x$ &$y(x)$&$y'(x)$   \\ 
\midrule

0.1&0.998334582651024 &-0.033283374960220\\
0.2&0.993353288961344 &-0.066267995319313\\
0.3&0.985100745872271 &-0.098660068556290\\
0.4&0.973650509840501 &-0.130175582648867\\
0.5&0.959103856956817 &-0.160544891813613\\
0.6&0.941588132070691 &-0.189516931926819\\
0.7&0.921254699087677 &-0.216862968455471\\
0.8&0.898276543103152 &-0.242379797978458\\
0.9&0.872845582616537 &-0.265892334576062\\
1.0&0.845169755493675 &-0.287255540026184\\
2.0&0.495936764048973 &-0.372832141746160\\
3.0&0.158857608676200 &-0.284252727750886\\
3.6&0.011090994555729 &-0.209392664698195\\

\bottomrule
\end{tabular}
\end{table}
\end{center}
%%%%%%%%%%%%%%%%%%%%%%%%%%%%%%%%%%%%%%%%%%

%%%%%%%%%%%%%%%%%%%%%%%%%%%%%%%%%%%%%%%%%%
\footnotesize 
\begin{center}
\begin{table}[ht!]
\caption{Obtained values of $y(x)$ and $y'(x)$ of  standard Lane-Emden equations for m = 2 by basis of RBs  with $N=75$ and iterations 15 }
\centering
%\footnotesize 
\begin{tabular}[ht!]{c|l|l}
\toprule
   $x$ &$y(x)$&$y'(x)$   \\ 
\midrule

0.1&0.99833499854614841738470254242797205 &-0.03326675387428229266020296311565521\\
0.2&0.99335990717838006759696432600612293 &-0.06613611499043401334084315714516733\\
0.3&0.98513394694678774324661706112162525 &-0.09822101034279809290634911546527703\\
0.4&0.97375411632745104999586632039810568 &-0.12915454582043665332495150875712279\\
0.5&0.95935271580338270088050287823943708 &-0.15859897547445287273609058975673498\\
0.6&0.94209403572565310558222573899059883 &-0.18625341132766885551729875654275060\\
0.7&0.92217034852259238973227792180747576 &-0.21185998653016673029184608587598711\\
0.8&0.89979737025891753987119186898118925 &-0.23520828220779765886143053232378281\\
0.9&0.87520937032283832664335099313246129 &-0.25613793250692707919375214927915510\\
1.0&0.84865411140824967691151067041935559 &-0.27453942454799399072475655485844822\\
2.0&0.52983642933948885122298186592152665 &-0.32634885813595725790315062497173500\\
3.0&0.24182408305234091675614257310953074 &-0.24062145844675797515285563922897989\\
4.0&0.04884014997594444649562664567725640 &-0.15040965958043957223901919928204625\\
4.3&0.00681094327420583009083737134824214 &-0.13039647888956858753190354417241168\\

\bottomrule
\end{tabular}
\end{table}
\end{center}
%%%%%%%%%%%%%%%%%%%%%%%%%%%%%%%%%%%%%%%%%%

%%%%%%%%%%%%%%%%%%%%%%%%%%%%%%%%%%%%%%%%%%
\footnotesize 
\begin{center}
\begin{table}[ht!]
\caption{Obtained values of $y(x)$ and $y'(x)$ of  standard Lane-Emden equations for m = 2.5 by basis of RBs  with $N=75$ and iterations 15 }
\centering
%\footnotesize 
\begin{tabular}[ht!]{c|l|l}
\toprule
   $x$ &$y(x)$&$y'(x)$   \\ 
\midrule

0.1&0.998335414189491 &-0.033250148555062\\
0.2&0.993366508668235 &-0.066004732702853\\
0.3&0.985166960607077 &-0.097785664864449\\
0.4&0.973856692696194 &-0.128148702313160\\
0.5&0.959597754464204 &-0.156697706048055\\
0.6&0.942588917282480 &-0.183095996800778\\
0.7&0.923059301998553 &-0.207074283925069\\
0.8&0.901261395554722 &-0.228434944738734\\
0.9&0.877463820286722 &-0.247052726803513\\
1.0&0.851944199128236 &-0.262872200779799\\
2.0&0.558372334987405 &-0.290313683599236\\
3.0&0.306675101717593 &-0.208571050779423\\
4.0&0.137680733022609 &-0.134053438395795\\
5.0&0.029019186649369 &-0.087473533084964\\
5.3&0.004259543533703 &-0.077863974396729\\

\bottomrule
\end{tabular}
\end{table}
\end{center}
%%%%%%%%%%%%%%%%%%%%%%%%%%%%%%%%%%%%%%%%%%

%%%%%%%%%%%%%%%%%%%%%%%%%%%%%%%%%%%%%%%%%%
\footnotesize 
\begin{center}
\begin{table}[ht!]
\caption{Obtained values of $y(x)$ and $y'(x)$ of  standard Lane-Emden equations for m = 3 by basis of RBs  with $N=75$ and iterations 15 }
\centering
%\footnotesize 
\begin{tabular}[ht!]{c|l|l}
\toprule
   $x$ &$y(x)$&$y'(x)$   \\ 
\midrule

0.1&0.99833582956916949360303219005188124 &-0.03323355906978612543528108656894298\\
0.2&0.99337309351037690525196177958205868 &-0.06587384697328259582789650194519630\\
0.3&0.98519978859187628465790325123284167 &-0.09735398607439270035647436495914414\\
0.4&0.97395825591973737173183189357312441 &-0.12715771812656874931089211804643185\\
0.5&0.95983906994485172355015029196763836 &-0.15483957691106171525826169820184137\\
0.6&0.94307317270163244519537716363760024 &-0.18003963338652410623289623950802367\\
0.7&0.92392283802402754599371328152955595 &-0.20249208705551806529037852452154313\\
0.8&0.90267208912835443207389240453497315 &-0.22202765241699836917057715071349470\\
0.9&0.87961716706036921959305769658330247 &-0.23857028971216663116635342349178773\\
1.0&0.85505756858862631144699757905789541 &-0.25212927977768245416231686995936214\\
2.0&0.58285051510965197279807688358867210 &-0.26149092569261197744719124891485520\\
3.0&0.35922650065961804953058851682521728 &-0.18404987913139678752701432127710201\\
4.0&0.20928161332783184750722991166749403 &-0.12016906000415030345933138030567754\\
5.0&0.11081983513962559885830939817050228 &-0.08012604337165627260209120489931835\\
6.0&0.04373798388970814005547798486698057 &-0.05604388226451708827307318053155939\\
6.8&0.00416778936545346001263963263378843 &-0.04364696951001550395471856003886447\\

\bottomrule
\end{tabular}
\end{table}
\end{center}
%%%%%%%%%%%%%%%%%%%%%%%%%%%%%%%%%%%%%%%%%%

%%%%%%%%%%%%%%%%%%%%%%%%%%%%%%%%%%%%%%%%%%
\footnotesize 
\begin{center}
\begin{table}[ht!]
\caption{Obtained values of $y(x)$ and $y'(x)$ of  standard Lane-Emden equations for m = 4 by basis of RBs  with $N=75$ and iterations 15 }
\centering
%\footnotesize 
\begin{tabular}[ht!]{c|l|l}
\toprule
   $x$ &$y(x)$&$y'(x)$   \\ 
\midrule

0.1&0.99833665953957353917 &-0.03320042731101602052\\
0.2&0.99338621353236887458 &-0.06561355430127865539\\
0.3&0.98526489445824457228 &-0.09650144694916813609\\
0.4&0.97415840895070184085 &-0.12521904232653407185\\
0.5&0.96031090234222125391 &-0.15124704523040264218\\
0.6&0.94401129085560210481 &-0.17421139290379387733\\
0.7&0.92557835269653368985 &-0.19388869549916036586\\
0.8&0.90534592383779093911 &-0.21019908106443456806\\
0.9&0.88364932397603694257 &-0.22318930318706216396\\
1.0&0.86081381220831175185 &-0.23300964460615518736\\
2.0&0.62294077167068319754 &-0.21815323531073192916\\
3.0&0.44005069158766127850 &-0.14895436785082222650\\
4.0&0.31804242903566436744 &-0.09886802020831413214\\
5.0&0.23592273104248679739 &-0.06788810347440624083\\
6.0&0.17838426534298279218 &-0.04865643577466167176\\
7.0&0.13635230535983164961 &-0.03626805424834208635\\
8.0&0.10450408207160914867 &-0.02795075318477840998\\
9.0&0.07961946745395432400 &-0.02214833117831084820\\
10 &0.05967274158948932881 &-0.01796142023434323612\\
11 &0.04334009538193507922 &-0.01485063006054293705\\
12 &0.02972593235798682964 &-0.01248033393137584648\\
13 &0.01820540390617142867 &-0.01063445527740952134\\
14 &0.00833052669542489543 &-0.00916953946501606750\\
14.9&0.00057641886621354664 &-0.00809526559361695336\\

\bottomrule
\end{tabular}
\end{table}
\end{center}
%%%%%%%%%%%%%%%%%%%%%%%%%%%%%%%%%%%%%%%%%%

\end{document}